\documentclass[11pt]{article}
\usepackage{amsmath}
\usepackage{amssymb}
\usepackage{amsfonts}
\usepackage{mathrsfs}
\usepackage{shortvrb,psfrag}
\usepackage{enumerate}
\usepackage{color}
\usepackage[dvips]{graphicx}

\let\pa=\partial

\def\ba{\begin{eqnarray}}
\def\ea{\end{eqnarray}}
\def\mv{{\mathbf{v}}}
\def\mm{\mathbf{m}}
\def\mq{\mathbf{Q}}
\def\mn{{\mathbf{n}}}
\def\md{{\mathbf{D}}}
\def\momega{{\mathbf{\Omega}}}

\def\mh{{\mathbf{H}}}
\def\mM{{\mathbf{M}}}
\def\mI{\mathbf{I}}

\def\cJ{{\cal J}}

\def\R{\Bbb R}

\def\no{\noindent}

\def\endproof{\hphantom{MM}\hfill\llap{$\square$}\goodbreak}

\newcommand{\beq}{\begin{equation}}
\newcommand{\eeq}{\end{equation}}
\newcommand{\ben}{\begin{eqnarray}}
\newcommand{\een}{\end{eqnarray}}
\newcommand{\beno}{\begin{eqnarray*}}
\newcommand{\eeno}{\end{eqnarray*}}

\newtheorem{Theorem}{Theorem}[section]
\newtheorem{Definition}[Theorem]{Definition}

\newtheorem{Lemma}[Theorem]{Lemma}

\newtheorem{Remark}[Theorem]{Remark}

\begin{document}

\title{From the Q-tensor flow for the liquid crystal to the Harmonic map flow }

\author{Meng Wang$^\dag$, Wendong Wang$^\ddag$ and Zhifei Zhang$^\sharp$
\\[2mm]
{\small $ ^\dag$  Department of Mathematics, Zhejiang University, Hangzhou
310027, P. R. China}\\
{\small E-mail:  mathdreamcn@zju.edu.cn}\\[2mm]
{\small $ ^\ddag$ School of Mathematical Sciences, Dalian University of Technology, Dalian, 116024,  China}\\
{\small E-mail: wendong@dlut.edu.cn}\\[2mm]
{\small $ ^\sharp$ School of  Mathematical Sciences and LMAM, Peking University, Beijing 100871, China}\\
{\small E-mail: zfzhang@math.pku.edu.cn}}


\date{\today}
\maketitle

\begin{abstract}
In this paper, we consider the  solutions of the relaxed Q-tensor flow in $\R^3$ with small parameter $\epsilon$. Firstly, we show that the limiting map is the so called harmonic map flow; Secondly, we also present a new proof for the global existence of weak solution for the harmonic map flow in three dimensions as in \cite{struwe88} and \cite{keller}, where Ginzburg-Landau approximation approach was used.
\end{abstract}

\setcounter{equation}{0}
\section{Introduction}
Liquid crystals are a state of matters that have properties between those of a conventional liquid and those of a solid  crystal. One of the most common liquid crystal phases is the nematic. The nematic liquid crystals are composed of rod-like molecules with the long axes of neighboring molecules approximately aligned to  one another. There are three different kinds of theories to model the nematic liquid crystals: Doi-Onsager theory, Landau-de Gennes theory and Ericksen-Leslie theory. The first is the molecular kinetic theory, and the later two are the continuum theory. In the spirit of Hilbert sixth problem, it is very important to explore the relationship between these theories.

Ball-Majumdar\cite{bm} define a Landau-de Gennes type energy functional in terms of the mean-field Maier-Saupe energy. Majundar-Zarnescu\cite{mz} consider the Oseen-Frank limit of the static  Q-tensor model. Their results show that the predictions of the Oseen-Frank theory and the Landau-De Gennes theory agree away from the singularities of the limiting Oseen-Frank global minimizer.

In \cite{kd, ez}, Kuzzu-Doi and E-Zhang formally derive the Ericksen-Leslie equation from the Doi-Onsager equations by taking small Deborah number limit. In \cite{wzz1,wzz2}, Wang-Wang-Zhang present a rigorous derivation from Doi-Onsager theory and Landau-de Gennes theory. In \cite{hlwz}, a systematical approach was proposed to derive the continuum theory from the molecular kinetic theory in both static and dynamic case.

Different with the above results on the static or local case, the goal of this work is to investigate the global convergence problem of  the solutions $Q_{\epsilon}$ of the Q-tensor flow in  $\R^3$. We will show that $Q_{\epsilon}$  convergences weakly to the weak solution of the harmonic map flow.

\subsection{The relaxed  Q-tensor flow}
In Landau-de Gennes theory, the state of the nematic liquid crystals is described by the macroscopic Q-tensor order parameter, which is a symmetric, traceless $3\times 3$ matrix. Physically, it can be interpreted as the second-order moment of the orientational distribution  function $f$, that is
\begin{eqnarray}
Q=\int_{S^2}(\mm\mm-\frac{1}{3}Id)f d\mm.\nonumber
\end{eqnarray}
 When $Q=0$, the nematic liquid crystal is said to be {\bf isotropic}. When $Q$ has two equal non-zero eigenvalues, it is said to be {\bf uniaxial} and $Q$ can be written as
 \begin{eqnarray}
 Q=s(\mn\mn-\frac{1}{3}Id),~~\mn\in S^2.\nonumber
  \end{eqnarray}
 When $Q$ has three distinct eigenvalues, it is said to be {\bf biaxial} and $Q$ can be written as
 \begin{eqnarray}
 Q=s(\mn\mn-\frac{1}{3}Id)+\lambda(\mn'\mn'-\frac{1}{3}Id),~\mn,\mn'\in S^2,~~\mn\cdot\mn'=0.\nonumber
  \end{eqnarray}
The general Landau-de Gennes energy functional takes the form
\beno
{\cal{F}_{LG}}(Q,\nabla Q)&=&\int_{\R^3}\Big\{\underbrace{-\frac{a}{2}tr\mq^2-\frac{b}{3}tr\mq^3+\frac{c}{4}(tr\mq^2)^2}_{f_{B}:\mbox{bulk energy}}\nonumber\\
&&+\frac{1}{2}\underbrace{\left(L_1|\nabla\mq|^2+L_2Q_{ij,j}Q_{ik,k}+L_3Q_{ij,k}Q_{ik,j}+L_4Q_{ij}Q_{kl,i}Q_{kl,j}\right)}_{f_{E}:\mbox{elastic energy}}\Big\}dx,
\eeno
here $a,b,c$ are material-dependent and temperature-dependent nonnegative constant and $L_i(i=1,2,3,4)$ are material dependent elastic constants. We refer to \cite{g,mn} for more details.

There are several dynamic Q-tensor models to describe the flow of the nematic liquid crystal, which are either derived from the molecular kinetic theory for the rigid rods by various closure approximation such as \cite{fcl,fls,hlwz}, or directly derived by variational method such as Beris-Edwards model \cite{be} and Qian-Sheng's model \cite{qs}.

In \cite{wzz2}, the authors consider the following Beris-Edwards model
\begin{align} \label{eq:MHD}(\rm Q)\,\, \left\{
\begin{aligned}
&\frac{\partial\mv}{\pa t}+\mv\cdot\nabla\mv=-\nabla p+\nabla\cdot(\sigma^s+\sigma^a+\sigma^d),\\
&\nabla\cdot\mv=0,\\
&\frac{\partial Q}{\pa t}+\mv\cdot\nabla Q+Q\cdot\Omega-\Omega\cdot Q=\frac{1}{\Gamma}H+S_{Q}(D).
\end{aligned}
\right. \end{align}

Here $\Gamma$ is a collective rotational diffusion constant, and
$$\md=\frac{1}{2}\left(\nabla\mv+(\nabla\mv)^T\right),\quad \momega=\frac{1}{2}(\nabla\mv-(\nabla\mv)^T).$$
Moreover, $\sigma^s$, $\sigma^a$ and $\sigma^d$ are symmetric viscous stress, antisymmetric viscous stress and distortion stress respectively defined by
\beno
\sigma^s=\eta\md-S_{\mq}(\mh),\quad \sigma^a=\mq\cdot\mh-\mh\cdot\mq, \quad \sigma^{d}=-\frac{\partial{\cal{F}_{LG}}}{\partial Q_{kl,j}} Q_{kl,i},
\eeno
where $\eta>0$ is the viscous coefficient,
 $\mh$ is the molecular field given by
\ben
\mh(\mq)=-\frac{\delta{\cal{F}_{LG}}}{\delta \mq},
\een
and
$S_{\mq}(\mM)$ is defined by
$$
S_{\mq}(\mM)=\xi\left(\mM\cdot(\mq+\frac{1}{d}\mI)+(\mq+\frac{1}{d}\mI)\cdot\mM-2(\mq+\frac{1}{d}\mI)\mq:\mM\right)
$$
for symmetric and traceless matrix $M$, where $\xi$ is a constant depending on the molecular details of a given liquid crystals. The well-posedness results of the Q-tensor model are studied in \cite{pz1,pz2}.

Wang-Zhang-Zhang \cite{wzz2} justify the limit  from Beris-Edwards system with a small parameter $\epsilon$ to the Ericksen-Leslie system before the first singularity time of the limit system. The limit behavior of the solution after the singularity remains unknown.  In this paper, we are interested in the global convergence from the Q-tensor flow to the harmonic map flow in $\R^3$.  Let us begin with the simplest form $L_2=L_3=L_4=0$, i.e.,
\ben
\frac{\partial\mq}{\partial t}=-\frac{1}{\epsilon\Gamma}\mathcal{J}(Q)+\frac{L_1}{\Gamma}\Delta Q\label{equation:q},
\een
where
\beno
{\cJ}(\mq):=\frac{\delta f_{B}(\mq)}{\delta\mq}=-a\mq-b\mq^2+c|\mq|^2\mq+\frac{1}{3}b|\mq|^2\mI.
\eeno

Let  $\vec{b}\in  S^2$ is a constant vector, $\mn_0:\R^3\rightarrow S^2$ such that $
\mn_0-\vec{b}\in H^{s+1}(\R^3)(s>0)
$. Moreover, $Q_0(x)=s_{+}(\mn_0(x)\otimes\mn_0(x)-\frac{Id}{3})$. We consider the following relaxed Q-tensor equations with a small parameter $\epsilon$:
\begin{align} \label{qtensorflow}(\rm Q_{\epsilon})\,\, \left\{
\begin{aligned}
&\frac{\partial Q_{\epsilon}}{\partial t}=\frac{(a-c|Q_{\epsilon}|^2)Q_{\epsilon}+bQ_{\epsilon}^2-b\frac{|Q_{\epsilon}|^2}{3}Id}{\epsilon\Gamma}+\frac{L_1\Delta Q_{\epsilon}}{\Gamma},\\
&Q_{\epsilon}(\cdot,0)=Q_{0}(x),
\end{aligned}
\right. \end{align}
which has a unique strong solution $Q_{\epsilon}(t,x)$ satisfying
$Q_{\epsilon}(t,x)-Q_0(x)\in L_{loc}^{\infty}([0,\infty),\dot{H}^{s+1}(\R^3))$.
We will study the global convergence of $Q_{\epsilon}$ as $\epsilon$
tends to zero.

\subsection{Main result}

The initial data $Q_0$ of the Q-tensor flow equations (\ref{qtensorflow}) lies in a special space, which contains the minimizers of the bulk energy $f_{B}(Q)$. To begin this,  we introduce some notations and known results.

Let $M^{sym}_{3\times 3}$ denote the set of real $3\times 3$ symmetric matrices and      $\mathcal{Q}_0\subset M^{3\times 3}$ denote the space of $Q$-tensors defined by
\begin{eqnarray}
\mathcal{Q}_{0}:=\big\{Q\in M_{3\times 3}^{sym}, Q_{ii}=0\big\},\nonumber
\end{eqnarray}
where we have used the Einstein summation convention.
The matrix norm is defined as
\begin{eqnarray}
|Q|:=\sqrt{tr Q^2}=\sqrt{Q_{ij}Q_{ij}}.\nonumber
\end{eqnarray}
We also write
\begin{eqnarray}
|\nabla Q|^2=\partial_{\alpha}Q_{ij}\partial_{\alpha}Q_{ij}.\nonumber
\end{eqnarray}
The bulk energy density can be written as
\begin{eqnarray}\label{eq:f B formula}
f_{B}(Q)=-\frac{a}{2}|Q|^2-\frac{b}{3}tr(Q^3)+\frac{c}{4}|Q|^4.
\end{eqnarray}
One can verify that $f_{B}$ is bounded from below (for example, see \cite[Proposition 8]{mz}), thus $f_{B}$ has the corresponding non-negative bulk energy density $\tilde{f}_{B}$ defined by
\begin{eqnarray}\label{eq:tilde f B formula}
\tilde{f}_{B}(Q)=f_{B}(Q)-\min_{Q\in\mathcal{Q}_0}f_{B}(Q).
\end{eqnarray}
In \cite[Proposition 8]{mz}, it was proved that $\tilde{f}_{B}$ attains its minimum on the uniaxial $Q-$tensors with constant order parameter $s_{+}$  as shown below
\begin{eqnarray}\label{eq:Q0 in N}
\tilde{f}_{B}(Q)=0~\Leftrightarrow ~Q\in \mathcal{N}~~\mbox{where}\nonumber\\
\mathcal{N}=\left\{Q\in\mathcal{Q}_0,~Q=s_{+}\left(\mn\otimes\mn-\frac{1}{3}Id\right),\mn\in S^2\right\},\label{mathcaln}
\end{eqnarray}
with
\begin{eqnarray}\label{eq:s+}
s_{+}=\frac{b+\sqrt{b^2+24ac}}{4c}.
\end{eqnarray}
For a matrix $Q\in \mathcal{N}$, we use $T_{Q}\mathcal{N}$ to denote the tangent space to $\mathcal{N}$ at $Q$ in $\mathcal{Q}_0$, and $\left(T_{Q}\mathcal{N}\right)^{\bot}_{\mathcal{Q}_0}$ to denote the orthogonal complement of $T_{Q}\mathcal{N}$ to $\mathcal{Q}_0$.

Let $\{A,B\}=AB+BA$ for $A,B\in M_{3\times 3}^{sym}.$ It was described in \cite[Lemma 2]{nz} that
\ben
T_{Q}\mathcal{N}&=&\left\{\dot{Q}\in M_{3\times 3}^{sym}: \frac{1}{3}s_{+}\dot{Q}=\{\dot{Q},Q\}\right\},\nonumber\\
&=&\big\{\mn\otimes\dot{\mn}+\dot{\mn}\otimes\mn:\dot{\mn},\dot{\mn}\in T_{\mn}{S^2}\big\},
\een
and
\ben\label{bot}
(T_{Q}\mathcal{N})^{\bot}_{\mathcal{Q}_0}&=&\left\{Q^{\bot}\in\mathcal{Q}_0:Q^{\bot}Q=QQ^{\bot}\right\}.
\een
We will specifically describe the orthogonal basis for  the tangent and normal space in the following Lemma \ref{lem:normalspace}.  It is obvious from (\ref{bot}) that $\mathcal{J}(Q)\in (T_{Q}\mathcal{N})^{\bot}_{\mathcal{Q}_0}$ for $Q\in \mathcal{N}$, and we will show in Lemma \ref{lem:the equivalent distance} that  for the approximating $Q_{\epsilon}$ near $\mathcal{N}$, we still have
\beno
\mathcal{J}(Q_{\epsilon})\in (T_{\pi_{\mathcal{N}}(Q_{\epsilon})}\mathcal{N})_{\mathcal{Q}_0}^{\bot},\eeno
where $\pi_{\mathcal{N}}$ denotes the projection operator on $\mathcal{N}$.

Let $z=(x,t)$ denote points in $\R^3\times \R_{+}$. For $z_0=(x_0,t_0)$, $R>0$, let
\beno
P_{R}(z_0)=\{z=(x,t)| |x-x_0|<R, |t-t_0|<R^2\},\\
S_{R}(z_0)=\{z=(x,t)|t=t_0-R^2\},\\
T_{R}(z_0)=\{z=(x,t)|t_0-4R^2<t<t_0-R^2\}.
\eeno
Denote the scaled fundamental solution to the heat equation
\begin{eqnarray}
G_{z_0}(z)=\tilde{G}_{(\frac{\Gamma}{L_1} x_0,\frac{\Gamma}{L_1} t_0)}\left(\frac{\Gamma}{L_1}x,\frac{\Gamma}{L_1}t\right),\nonumber
\end{eqnarray}
where
\begin{eqnarray}
\tilde{G}_{z_0}(z)=\frac{1}{(4\pi(t_0-t))^{3/2}}e^{-\frac{|x-x_0|^2}{4(t_0-t)}}, t<t_0.\nonumber
\end{eqnarray}
Then we have
\beno
\nabla G_{z_0}=\frac{\Gamma}{2L_1}\cdot\frac{x-x_0}{t-t_0} G
\eeno
Also we write $P_{R}(0)=P_{R}$, $T_{r}(0)=T_{r}$, and $G_{0}(z)=G(z)$.

Similar to the harmonic map flow, the limiting Q-tensor flow takes as follows
\ben\label{eq:limit Q tensor flow}
\partial_t Q-\frac{L_1}{\Gamma}\Delta Q+\lambda(x,t)\gamma_{\mathcal{N}}(Q)=0,\quad Q(x,t)|_{t=0}=Q_0(x),\een
where $\lambda(x,t)$ is a function of $L^2_{loc}$, and  $\gamma_{\mathcal{N}}(Q)$ is unit normal vector to $(T_{Q}\mathcal{N})^{\bot}_{\mathcal{Q}_0}$ at $Q$.

\begin{Definition}
A $Q-$tensor $Q(x,t):\R^3\times\R_{+}\to \mathcal{N}$ is called a weak solution to (\ref{eq:limit Q tensor flow}) if $Q(x,t)|_{t=0}=Q_0(x)$ a.e., $\partial_tQ,\nabla Q\in L_{loc}^{2}(\R^3\times \R_{+})$, and there holds
\begin{eqnarray}
\int_{\R^3}\int_{\R_{+}}(\partial_t Q:\phi+\frac{L_1}{\Gamma}\nabla Q:\nabla\phi+\lambda\gamma_{\mathcal{N}}(Q):\phi)dxdt=0,\nonumber
\end{eqnarray}
for all $\phi\in C_{0}^{\infty}(\R^3\times\R_{+},\R^{3\times 3})$.
\end{Definition}

Our main theorem is stated as follows.

\begin{Theorem}\label{thm:main}
Let $Q_{\epsilon}$ satisfy the equations of the relaxed Q-tensor flow equations (\ref{qtensorflow}) with the data $Q_0(x)=s_{+}\left(\mn_0\otimes\mn_0-\frac{1}{3}Id\right)\in \mathcal{N}$ as in (\ref{eq:Q0 in N}). Then
\item 1. there exists a subsequence of $\epsilon$ (also denoted $\epsilon$) such that
\beno
Q_{\epsilon}\rightharpoonup Q=s_{+}\left(\mn\otimes\mn-\frac{1}{3}Id\right)\in \mathcal{N},
\eeno
and $Q$ solves weakly the equations  (\ref{eq:limit Q tensor flow}).\\
\item 2. the director field $\mn$ weakly solves
\begin{eqnarray}\label{heatflow}
\partial_t\mn-\Delta\mn=|\nabla\mn|^2\mn,\quad \mn(x,t)|_{t=0}=\mn_0(x).
\end{eqnarray}
\item 3. $\mn$ is regular on a dense open set $\Omega_0\subset \R^{3}\times\R_{+}$, whose complement  $\Sigma$ has locally finite 3-dimensional Hausdorff-measure (with respect to the parabolic metric).
\end{Theorem}

\begin{Remark}
Compared with \cite{wzz1,wzz2}, the above theorem makes it reasonable that the global weak convergence from the Q-tensor flow to the simplified Oseen-Frank map flow (i.e, harmonic map flow). Different from Ginzburg-Landau  approximation used in \cite{struwe88} and \cite{keller}, we consider the Q-tensor approximation, and the difficulty is that the properties of the limit manifold $\mathcal{N}$  is unclear as stated in \cite{nz}. In the next section, we give a careful study for the geometry of $\mathcal{N}$ ( see Lemma \ref{lem:normalspace}).
\end{Remark}

\setcounter{equation}{0}
\section{Technical lemmas and interior regularity estimates}
In this section, we will introduce the properties of Q-tensor matrix,  the tangent space, the normal space and the equivalence of the bulk energy. Using these estimates and exploring monotonicity inequalities as in \cite{struwe88}, we can obtain the interior regularity criteria of Q-tensor equations.

First of all, for the matrix of $3\times 3$, we have the following properties.
\begin{Lemma}\label{lem:AB}
 Let $A,B$ be matrices of $3\times 3$.\\
(i) If $A$ is symmetric, then
$$A:B=A:\tilde{B},\quad \tilde{B}=\frac{B+B^{T}}{2}$$
where
$\tilde{B}$ is the symmetrization for $B$.\\
(ii)If $A$ is antisymmetric, then
$$A:B=A:\bar{B},\quad \bar{B}=\frac{B-B^{T}}{2}$$
where
$\bar{B}$ is the antisymmetrization for $B$.\\
(iii) If $A$ is symmetric and $B$ is antisymmetric, then
$$ A:B=0.$$
\end{Lemma}

For $Q\in \mathcal{N}$, it is easy to verify that the orthogonal basis of  $T_{Q}\mathcal{N}$ and $(T_{Q}\mathcal{N})^{\bot}_{\mathcal{Q}_0}$ is as follows.

\begin{Lemma}\label{lem:normalspace}
Let $Q=s_{+}\left(\mn_3\otimes\mn_3-\frac{1}{3}Id\right)\in \mathcal{N}$, and $\mn_1,\mn_2$ be unit perpendicular vectors in $ V_{\mn_3}=\{\mn^{\bot}\in \R^3: \mn^{\bot}\cdot\mn_3=0\}$. Then it holds that
\begin{itemize}
\item[1.]
\ben\label{eq:tangent space}
T_{Q}\mathcal{N}&=&Span\left\{\frac{1}{\sqrt{2}}(\mn_3\otimes\mn_2+\mn_2\otimes\mn_3),\frac{1}{\sqrt{2}}(\mn_3\otimes\mn_1+\mn_1\otimes\mn_3)\right\},
\een
\item[2.]
\ben\label{eq:normal space}
\left(T_{Q}\mathcal{N}\right)^{\bot}_{\mathcal{Q}_0}&=&Span\left\{\frac{1}{\sqrt{2}}(\mn_2\otimes\mn_1+\mn_1\otimes\mn_2), \frac{1}{\sqrt{2}}(\mn_1\otimes\mn_1-\mn_2\otimes\mn_2),\right.\nonumber\\
&&~~\left.\sqrt{6}\left(\frac{1}{2}\mn_1\otimes\mn_1+\frac{1}{2}\mn_2\otimes\mn_2-\frac{Id}{3}\right)\right\}
\een
\item[3.]Moreover,
$
{\mathcal Q}_{0}=T_{Q}\mathcal{N}\oplus (T_{Q}\mathcal{N})_{\mathcal{Q}_0}^{\bot}.
$
\end{itemize}
\end{Lemma}

\no{\it Proof:} The tangent space at $s_{+}(\mn_3\otimes\mn_3-\frac{1}{3}Id)$ of $\mathcal{N}$ (\ref{eq:tangent space}) is a direct result from  \cite[(2.3)]{wzz2}. The others can be deduced by direct computations.\endproof

\begin{Lemma}\label{lem:Gamma Q in N} For $Q\in \mathcal{Q}_0$,
there exists $\epsilon_0>0$ such that if $dist(Q,\mathcal{N})<\epsilon_0$, then
\begin{eqnarray}\label{eq:Gamma Q in N}
{\mathcal{J}}(Q)\in (T_{\pi_{\mathcal{N}}(Q)}(\mathcal{N}))_{\mathcal{Q}_0}^{\bot}.
\end{eqnarray}
\end{Lemma}
{\it Proof:}
Denote the eigenvectors of $Q(x,t)$ by $\mn_1(x,t),\mn_2(x,t),\mn_{3}(x,t)$ corresponding to its eigenvalues $\lambda_{1}(x,t)$, $\lambda_{2}(x,t)$, $\lambda_{3}(x,t)=-\lambda_1(x,t)-\lambda_2(x,t)$. Then we have
\ben\label{eq:Q identity}
Q(x,t)=\lambda_{1}\mn_1\otimes\mn_1+\lambda_{2}\mn_2\otimes\mn_2+\lambda_{3}\mn_3\otimes\mn_3,
\een
especially,
\ben\label{eq:I identity}
I=\mn_1\otimes\mn_1+\mn_2\otimes\mn_2+\mn_3\otimes\mn_3.
\een

Choose $\epsilon_0$ small enough such that $dist(Q,\mathcal{N})<\epsilon_0$, then
\begin{eqnarray}\label{eq:Q N distance}
dist(Q,\mathcal{N})^2=\left(\lambda_1+\frac{s_{+}}{3}\right)^2+\left(\lambda_2+\frac{s_+}{3}\right)^2+\left(\lambda_1+\lambda_2+2\frac{s_{+}}{3}\right)^2,
\end{eqnarray}
furthermore,
$$dist(Q,\mathcal{N})^2=|Q-\pi_{\mathcal{N}}(Q)|^2,\quad \pi_{\mathcal{N}}(Q)=s_{+}(\mn_3\otimes\mn_3-\frac{Id}{3})$$
where $\pi_{\mathcal{N}}(Q)$ is unique and depends continuously on $Q$. See \cite[Lemma 8]{nz} for more details.

It is easy to see the projection of $Q$ on $T_{s_{+}(\mn_3\otimes\mn_3-\frac{1}{3}Id)}\mathcal{N}$ is $0$. Using (\ref{eq:I identity}), we have
\beno
Q&=&\lambda_1\mn_1\otimes\mn_1+\lambda_2\mn_1\otimes\mn_2-(\lambda_1+\lambda_2)\mn_3\otimes\mn_3\nonumber\\
&=&(2\lambda_1+\lambda_2)\left(\mn_1\otimes\mn_1-\frac{Id}{3}\right)+(\lambda_1+2\lambda_2)\left(\mn_2\otimes\mn_2-\frac{Id}{3}\right),
\eeno
\beno
Q^2&=&\lambda_1^2\mn_1\otimes\mn_1+\lambda_2^2\mn_2\otimes\mn_2+(\lambda_1+\lambda_2)^2\mn_3\otimes\mn_3\nonumber\\
&=&-\lambda_2(2\lambda_1+\lambda_2)\left(\mn_1\otimes\mn_1-\frac{Id}{3}\right)-\lambda_1(\lambda_1+2\lambda_2)\left(\mn_2\otimes\mn_2-\frac{Id}{3}\right)\nonumber\\
&&+\frac{(\lambda_{1}^{2}+\lambda_{2}^2+\lambda_{3}^{2})}{3}Id,
\eeno
which together with (\ref{eq:Q identity}) yields that
\beno
Q^2-\frac{|Q|^2}{3}Id=-\lambda_2(2\lambda_1+\lambda_2)\left(\mn_1\otimes\mn_1-\frac{Id}{3}\right)-\lambda_1(\lambda_1+2\lambda_2)\left(\mn_2\otimes\mn_2-\frac{Id}{3}\right).
\eeno
Hence, we get
\beno
{\mathcal{J}}(Q)&=&(a-c|Q|^2)Q+b\left(Q^2-\frac{|Q|^2}{3}Id\right)\\
&=&(a-c|Q|^2-b\lambda_2)(2\lambda_1+\lambda_2)(\mn_1\otimes\mn_1-\frac{Id}{3})+(a-c|Q|^2-b\lambda_1)(\lambda_1+2\lambda_2)(\mn_2\otimes\mn_2-\frac{Id}{3})\nonumber\\
&=&\frac{1}{\sqrt{2}}\left(a-c|Q|^2-b\lambda_3\right)(\lambda_1-\lambda_2)e_2\nonumber\\
&&+\frac{1}{\sqrt{6}}\left((a-c|Q|^2-b\lambda_2)(2\lambda_1+\lambda_2) +(a-c|Q|^2-b\lambda_1)(\lambda_1+2\lambda_2)\right)e_3,
\eeno
where we have used the orthogonal basis $e_2$ and $e_3$ in the normal space $(T_{\pi_{\mathcal{N}}(Q)}(\mathcal{N}))_{\mathcal{Q}_0}^{\bot}$ (see (\ref{eq:normal space})), and
\beno
e_2=\frac{1}{\sqrt{2}}(\mn_1\otimes\mn_1-\mn_2\otimes\mn_2),\quad e_3=\sqrt{6}(\frac{1}{2}\mn_1\otimes\mn_1+\frac{1}{2}\mn_2\otimes\mn_2-\frac{Id}{3}).\eeno
Thus,  (\ref{eq:Gamma Q in N}) is an immediate result. \endproof
\medskip

In fact, the nonnegative bulk energy $\tilde{f}_{B}(Q)$ is equivalent to the distance from $Q$ to $\mathcal{N}$, which is stated as the following Lemma.
\begin{Lemma}\label{lem:the equivalent distance}
There exists $\epsilon_0>0$ such that if $dist(Q,\mathcal{N})<\epsilon_0$, then
\begin{eqnarray}\label{eq:the equivalent distance}
\frac{1}{C}dist(Q,\mathcal{N})^2\le \tilde{f}_{B}(Q)\le C(dist(Q,\mathcal{N}))^2,
\end{eqnarray}
where $C$ is independent of $Q$, but depends on $a,b,c$.
\end{Lemma}
Proof: Assume that the eigenvalues of $Q$ are $x,y,-x-y$. If $dist(Q,\mathcal{N})<\epsilon_0$ and $\epsilon_0$ is small enough, similar to (\ref{eq:Q N distance}), we have
\begin{eqnarray}\label{eq:Q N distance2}
dist(Q,\mathcal{N})^2=(x+\frac{s_+}{3})^2+(y+\frac{s_+}{3})^2+(x+y+2\frac{s_+}{3})^2\le \epsilon_1.
\end{eqnarray}
Let
\ben\label{eq:G}
(x+\frac{s_{+}}{3})^2+(y+\frac{s_{+}}{3})^2+(x+y+\frac{2s_{+}}{3})^2
\triangleq G(x,y).
\een

On the other hand, for the nonnegative bulk energy, we have
\ben\label{eq:fb H}
\tilde{f}_{B}(Q)
&=&-\frac{a}{2}|Q|^2-\frac{b}{3}tr Q^3+\frac{c}{4}|Q|^4-\min_{Q\in\mathcal{Q}_0}{f}_{B}(Q)\nonumber\\
&=&-a(x^2+y^2+xy)+b(x^2y+xy^2)+c(x^2+y^2+xy)^2-\min_{Q\in\mathcal{Q}_0}{f}_{B}(Q)\nonumber\\
&\triangleq&H(x,y),
\een
where $H:\R^2\to \R$  is the 2-dimensional function as $G(x,y)$. Note that $H(x,y)=0$ only at three pairs $(x,y)$ namely $(-\frac{s_{+}}{3},-\frac{s_{+}}{3})$, $(2\frac{s_{+}}{3},-\frac{s_{+}}{3})$ and $(-\frac{s_{+}}{3},2\frac{s_{+}}{3})$, c.f.\cite[Lemma 5]{mz}. This fact together with $(\ref{eq:G})$ and $(\ref{eq:fb H})$ gives
\ben\label{eq:H G analysis}
H\left(-\frac{s_{+}}{3},-\frac{s_{+}}{3}\right)=\frac{\partial H}{\partial x}\left(-\frac{s_{+}}{3},-\frac{s_{+}}{3}\right)=\frac{\partial H}{\partial y}\left(-\frac{s_{+}}{3},-\frac{s_{+}}{3}\right)=0,\nonumber\\
G\left(-\frac{s_{+}}{3},-\frac{s_{+}}{3}\right)=\frac{\partial G}{\partial x}\left(-\frac{s_{+}}{3},-\frac{s_{+}}{3}\right)=\frac{\partial G}{\partial y}\left(-\frac{s_{+}}{3},-\frac{s_{+}}{3}\right)=0.
\een
Careful computations show that
\beno
\frac{\partial H}{\partial x}&=&(2x+y)[-a+by+2c(x^2+y^2+xy)],\\
\frac{\partial^2 H}{\partial x\partial y}&=&[-a+by+2c(x^2+y^2+xy)]+(2x+y)[b+2c(2y+x)],\\
\frac{\partial^2 H}{\partial x^2}&=&2[-a+by+2c(x^2+y^2+xy)]+2c(2x+y)^2.
\eeno
Noting $2cs_{+}^{2}-bs_{+}-3a=0$, we have
 \beno
\frac{\partial^2H}{\partial x^2}\left(-\frac{s_{+}}{3},-\frac{s_{+}}{3}\right)=\frac{\partial^2H}{\partial y^2}\left(-\frac{s_{+}}{3},-\frac{s_{+}}{3}\right)=2cs_{+}^2=bs_{+}+3a.
\eeno
\beno
\frac{\partial^2H}{\partial x\partial y}\left(-\frac{s_{+}}{3},-\frac{s_{+}}{3}\right)=-bs_{+}+2cs_{+}^{2}=3a,
\eeno
\beno
\frac{\partial^2 G}{\partial x^2}=4=\frac{\partial^2 G}{\partial y^2},~~\frac{\partial^2G}{\partial x\partial y}=2.
\eeno
Then, for $(x,y)\neq(-\frac{s_{+}}{3},-\frac{s_{+}}{3})$, by the above computations and (\ref{eq:H G analysis}), we have
\begin{eqnarray}\label{eq:G H}
\frac{H(x,y)}{G(x,y)}=\frac{H_1(x,y)+R_H(x,y)}{G_{1}(x,y)},
\end{eqnarray}
where
\begin{eqnarray}\label{eq:H1}
H_{1}(x,y)=(bs_{+}+3a)\left((x+\frac{s_{+}}{3})^2+(y+\frac{s_{+}}{3})^2\right)+6a\left(x+\frac{s_{+}}{3}\right)\left(y+\frac{s_{+}}{3}\right),
\end{eqnarray}
\begin{eqnarray}\label{eq:G1}
G_{1}(x,y)=4\left(x+\frac{s_{+}}{3}\right)^2+4\left(x+\frac{s_{+}}{3}\right)\left(y+\frac{s_{+}}{3}\right)+4\left(y+\frac{s_{+}}{3}\right)^2,
\end{eqnarray}
and $R_{H}$ is  the remainder in the Taylor expansions of $H(x,y)$ at $(-\frac{s_{+}}{3},-\frac{s_{+}}{3})$. Thus for $\epsilon_1$ sufficiently small in (\ref{eq:Q N distance2}), we get
\begin{eqnarray}\label{eq:RH}
|R_{H}(x,y)|\le \frac{bs_{+}}{2}\left((x+\frac{s_{+}}{3})^2+(y+\frac{s_{+}}{3})^2\right).
\end{eqnarray}
Summing up the inequalities (\ref{eq:G H})-(\ref{eq:RH}), we conclude that
$$
\frac{bs_{+}}{12}\leq \frac{H(x,y)}{G(x,y)}\leq (3a+\frac{3bs_{+}}{4}).
$$
The proof is completed. \endproof


Now we consider the evolution of the energy, and we will follow the same line as in \cite{struwe88}. First of all,  we define the energy density by
\begin{eqnarray}
e_{\epsilon}(Q,\nabla Q)=\frac{1}{\epsilon\Gamma}\tilde{f}_{B}(Q)+\frac{L_1}{2\Gamma}|\nabla Q|^2.
\end{eqnarray}

For the equations of (\ref{qtensorflow}) with initial data $Q_0\in \mathcal{N}$, there exists a global weak solution $Q_{\epsilon}(x,t)$ (denoted by $Q$ for simplicity), see \cite{pz2}. The solution $Q$ is regular indeed by usual energy estimates, and we have the following basic estimates.

\begin{Lemma}\label{lem:energy inequality}(Energy Inequality) Suppose that $Q(x,t)$ solves (\ref{qtensorflow}) with initial data $Q_0\in \mathcal{N}$, then it holds that
\begin{eqnarray}\sup_{t\ge 0}\left[\int_{0}^{t}\int_{\R^3}|\partial_t Q|^2dxdt+\int_{\R^3}e_{\epsilon}(Q(\cdot,t),\nabla Q(\cdot,t)) dx\right]\le \frac{L_1}{2\Gamma}\int_{\R^3}|\nabla Q_0|^2dx.
\end{eqnarray}
\end{Lemma}
Proof: Multiplying (\ref{qtensorflow}) by $Q_{t}$ and integration by parts yield that
\beno
\int_{\R^3}|\partial_t Q|^2=\frac{1}{\epsilon\Gamma}\int_{\R^3}-\frac{\delta \tilde{{f}}_{B}}{\delta Q}:\partial_t Q dx+\frac{L_1}{\Gamma}\int_{\R^3}\Delta Q:\partial_t Qdx.
\eeno
Noting that $\frac{\delta \tilde{{f}}_{B}}{\delta Q}:\partial_t Q=\partial_t  \tilde{f}_{B}( Q)$ and $\tilde{{f}}_{B}(Q_0)=0$, the lemma follows.\endproof
\medskip

The following parabolic maximal principle lemma is similar to Proposition 3 in \cite{mz}, where the elliptic case was considered. We omitted the proof.
\begin{Lemma}\label{lem:maximum principle}(Maximal Principle)Suppose that $Q(x,t)$ solves (\ref{qtensorflow}) with initial data $Q_0\in \mathcal{N}$, then it holds that
$$|Q|\le \sqrt{\frac{2}{3}}s_{+}.$$
\end{Lemma}

For the case of $L_2=L_3=L_4=0$, we also have the monotonicity properties of the level energy as in \cite{struwe88}.
\begin{Lemma}\label{lem:monotonicity} Suppose that $Q(x,t)$ solves (\ref{qtensorflow}) with initial data $Q_0\in \mathcal{N}$.
For any point $z_0=(x_0,t_0)\in \R^3\times \R_{+}$, the functions
\ben\label{eq:monotonicity}
\Phi(R,Q,\epsilon)=\frac{1}{\Gamma}R^2\int_{S_{R}(z_0)}\big[\frac{L_1}{2}|\nabla Q|^2+\frac{{\tilde f}_{B}(Q)}{\epsilon}\big]G_{z_0}dx,
\een
\ben\label{eq:monotonicity2}
\Psi(R,Q,\epsilon)=\frac{1}{\Gamma}\int_{T_{R}(z_0)}\big[\frac{L_1}{2}|\nabla Q|^2+\frac{{\tilde f}_{B}(Q)}{\epsilon}\big]G_{z_0}dxdt
\een
are non-decreasing for $0<R<\sqrt{t_0}/2$.
\end{Lemma}

\no{\it Proof:} We first note that $Q(x+x_0,t+t_0)$ satisfies the equations of (\ref{qtensorflow})  with initial data
$Q(\cdot,-t_0)=Q_{0}(x)$. Thus, we may assume that $z_0=(0,0)$. By scale invariance  $Q_R(x,t)=Q(Rx,R^2t)$ satisfying (\ref{qtensorflow}) with constant $\epsilon_{R}=\frac{\epsilon}{R^2}$, we have
\beno
\Phi(R,Q,\epsilon)&=&\frac{1}{\Gamma}R^2\int_{S_{R}}\big[\frac{L_1}{2}|\nabla Q|^2+\frac{{\tilde f}_{B}(Q)}{\epsilon}\big]Gdx\\
&=&\frac{1}{\Gamma}\int_{S_{1}}\big[\frac{L_1}{2}|\nabla Q_{R}|^2+\frac{{\tilde f}_{B}(Q_{R})}{\epsilon/R^2}\big]Gdx\triangleq\Phi(1,Q_{R},\epsilon/R^2).
\eeno
It suffices to consider the case of $R=1$. Direct computations and the equations (\ref{qtensorflow}) show that
\beno
&&\frac{d}{dR}\Phi(R,Q,\epsilon)|_{R=1}=\frac{d}{dR}\Phi(1,Q_R,\epsilon/R^2)|_{R=1}\nonumber\\
&=&\frac{1}{\Gamma}\int_{S_1}\left(L_1(-\Delta Q):(x\cdot\nabla Q+2t\partial_t Q)+\frac{1}{\epsilon}\frac{\delta \tilde{f}_{B}}{\delta Q}:(x\cdot\nabla Q+2t\partial_t Q)\right)G(x,-1)dx\nonumber\\
&&-\frac{1}{\Gamma} \int_{S_1}L_1(x\cdot\nabla Q+2t\partial_{t}Q):\nabla_k Q\nabla_k Gdx+\frac{2}{\epsilon\Gamma}\int_{S_1}\tilde{f}_{B}G(x,-1)dx\nonumber\\
&=&\int_{S_1}-\partial_tQ:(x\cdot\nabla Q+2t\partial_tQ)Gdx+\frac{1}{\Gamma}\int_{S_1}\frac{\Gamma G}{2t}(x\cdot\nabla Q):(x\cdot\nabla Q+2t\partial_tQ)dx\\
&&+\frac{2}{\epsilon\Gamma}\int_{S_1}\tilde{f}_{B}G(x,-1)dx\\
&=&\int_{S_1}\frac{1}{2|t|}\left(2t\partial_t Q+x\cdot\nabla Q\right)^2Gdx+\frac{2}{\epsilon\Gamma}\int_{S_1}\tilde{f}_{B}Gdx\ge 0,
\eeno
which implies the first inequality (\ref{eq:monotonicity}) for $R<\sqrt{t_0}$.

For $0<R<R_1<\sqrt{t_0}/2$, we consider the term $\Psi(R,Q,\epsilon)$ with $r'/r=R_1/R$, then
\beno
\Psi(R,Q,\epsilon)&=&\frac{1}{\Gamma}\int_{T_{R}}\big[\frac{L_1}{2}|\nabla Q|^2+\frac{{\tilde f}_{B}(Q)}{\epsilon}\big]G dxdt\\
&=&\frac{1}{\Gamma}\int_{-4R^2}^{-R^2}\int_{\R^3}\big[\frac{L_1}{2}|\nabla Q|^2+\frac{{\tilde f}_{B}(Q)}{\epsilon}\big]G dxdt\\
&=&2\int_{R}^{2R}r^{-1}\Phi(r,Q,\epsilon)dr\\
&=&2\int_{R}^{2R}\frac{\Phi(r,Q,\epsilon)}{\Phi(r',Q,\epsilon)}r'^{-1}\Phi(r',Q,\epsilon)dr'\leq \Psi(R_1,Q,\epsilon),
\eeno
where we used the monotonicity inequality (\ref{eq:monotonicity}). \endproof

\begin{Remark}\label{rem:monotonicity}
The above lemma indicates that the monotonic radius of $\Phi$ and $\Psi$ depends on $t_0$, which is reasonable since we have no definition for $t<0$. Similarly, if we consider the Q-tensor flow in $\R^3\times (-4R_0^2,R_0^2)$, then $\Phi(R,Q,\epsilon)$ is nondecreasing for $0<R<2R_0$ and
$\Psi(R,Q,\epsilon)$ is nondecreasing for $0<R<R_0.$
\end{Remark}

We have the following  Bochner-type inequality.

\begin{Lemma}\label{lem:bochner inequlity}
Suppose that $Q(x,t)$ solves (\ref{qtensorflow}) with initial data $Q_0\in \mathcal{N}$.
There exist  $\epsilon_0>0$ and a constant $C>0$, independent of $\epsilon$, such that
\begin{eqnarray}
(\partial_t-\Delta)e_{\epsilon}(Q,\nabla Q)(x,t)\le C e_{\epsilon}(Q,\nabla Q)^2(x,t),
\end{eqnarray}
provided that there exists a ball $B_{\rho(x)}(x)$ with $\rho(x)>0$ such that $\sup_{y\in B_{\rho(x)}(x)}dist(Q(y,t),\mathcal{N})<\epsilon_0$.
\end{Lemma}

\no{\it Proof:} Direct calculation shows that
\ben\label{eq:bochner1}
&&\left(\partial_t-\frac{L_1}{\Gamma}\Delta\right)\left(\frac{L_1}{2\Gamma}|\nabla Q|^2+\frac{\tilde{f}_{B}(Q)}{\epsilon\Gamma}\right)\nonumber\\
&=&\partial_{t}\left(\frac{L_1}{2\Gamma}|\nabla Q|^2+\frac{\tilde{f}_{B}}{\epsilon\Gamma}\right)-\frac{L_1^2}{2\Gamma^2}\Delta(|\nabla Q|^2)-\frac{L_1}{\epsilon\Gamma^2}\Delta(\tilde{f}_{B}(Q))\nonumber\\
&=&\frac{L_1}{\Gamma}\nabla Q:\partial_{t}\nabla Q+\frac{1}{\epsilon\Gamma}\frac{\delta \tilde{f}_{B}}{\delta Q}:\partial_tQ-\frac{L_1^2}{2\Gamma^2}\left(2\nabla Q:\Delta\nabla Q+2|\nabla^2Q|^2\right)\nonumber\\
&&-\frac{L_1}{\epsilon\Gamma^2}\frac{\delta\tilde{f}_{B}}{\delta Q}:\Delta Q-\frac{L_1}{\epsilon\Gamma^2}\nabla(\frac{\delta\tilde{f}_{B}}{\delta Q}):\nabla Q\nonumber\\
&=&-\frac{1}{\epsilon^2\Gamma^2}\left|\frac{\delta \tilde{f}_{B}}{\delta Q}\right|^2-\frac{2L_1}{\epsilon\Gamma^2}\nabla(\frac{\delta{\tilde{f}_{B}}}{\delta Q}):\nabla Q-\frac{L_1^{2}}{\Gamma^2}|\nabla^2 Q|^2.
\een
It suffices to estimate the second term of the above equality.

Denote the eigenvectors of $Q(x,t)$ by $\mn_1(x,t),\mn_2(x,t),\mn_{3}(x,t)$ corresponding to $\lambda_{1}(x,t),\lambda_{2}(x,t)$, $\lambda_{3}(x,t)=-\lambda_1(x,t)-\lambda_2(x,t)$. Then
\beno
\pi_\mathcal{N}(Q)=s_{+}\big(\mn_3(x,t)\otimes\mn_3(x,t)-\frac13 Id\big),
\eeno
which is a minimizer of the bulk energy $\tilde{f}_{B}$.
By the Taylor expansion of $\frac{\partial^2\tilde{f}_{B}}{\partial Q_{ij}\partial Q_{mn}}$ near $\pi_\mathcal{N}(Q)$, we get
\beno
&&\frac{\partial^2\tilde{f}_{B}}{\partial Q_{ij}\partial Q_{mn}}(Q(x,t))=\frac{\partial^2\tilde{f}_{B}}{\partial Q_{ij}\partial Q_{mn}}(\pi_\mathcal{N}(Q))\\
&&\quad \quad+\frac{\partial^3\tilde{f}_{B}}{\partial Q_{ij}\partial Q_{mn}\partial Q_{pq}}(\pi_\mathcal{N}(Q))(Q_{pq}(x,t)-\pi_\mathcal{N}(Q)_{pq})
+O(|Q_{pq}(x,t)-\pi_\mathcal{N}(Q)_{pq}|^2),
\eeno
where we have used the formula (\ref{eq:tilde f B formula}) of $\tilde{f}_{B}(Q)$ .

Using the convex property of  $\tilde{f}_{B}(Q)$ at  $\pi_\mathcal{N}(Q)$ and the maximum of $Q$ in Lemma \ref{lem:maximum principle}, we have
\beno
-\nabla(\frac{\delta {\tilde{f}}_{B}}{\delta Q}):\nabla Q&=&-\frac{\partial^2\tilde{f}_{B}}{\partial Q_{ij}\partial Q_{mn}}(Q(x,t))Q_{mn,k}(x,t)Q_{ij,k}(x,t)\\
&\le &C|Q(x,t)-\pi_\mathcal{N}(Q)|^2+C|\nabla Q|^4.
\eeno
Moreover, due to Lemma \ref{lem:the equivalent distance}, we get
\ben\label{eq:bochner2}
-\nabla(\frac{\delta {\tilde{f}}_{B}}{\delta Q}):\nabla Q &\le &C\tilde{f}_{B}(Q(x,t))+C|\nabla Q|^4.
\een
Combining (\ref{eq:bochner1}) and (\ref{eq:bochner2}), we obtain
\ben\label{eq:bochner3}
(\partial_t-\frac{\Gamma}{L_1}\Delta)e_{\epsilon}(Q)+\frac{1}{\epsilon^2\Gamma^2}\left|\frac{\delta \tilde{f}_{B}}{\delta Q}\right|^2\le Ce_{\epsilon}(Q)^2.
\een
The proof of the lemma is completed.\endproof

We consider the local uniform regularity property of the solution $Q_{\epsilon}$, which follows  from the monotonicity Lemma \ref{lem:monotonicity} and Schoen's trick, c.f. \cite[Theorem 5.1]{struwe88} or \cite[Theorem 2.2]{SchoenU}.

\begin{Lemma}\label{lem:interior regular} Suppose that $Q_{\epsilon}(x,t)$ solves (\ref{qtensorflow}) in $\R^3\times[-4R_0^2,R_0^2]$,
and there exist positive constants $\epsilon_0'$ and  $\epsilon_1$, such that when $\epsilon<\epsilon_0'$, for some $0<R<R_0$ the following inequality holds
\begin{eqnarray}
\Psi(R)=\Psi(R,Q_{\epsilon},\epsilon)=\int_{T_{R}}\left(\frac{L_1}{2\Gamma}|\nabla Q_{\epsilon}|^2+\frac{\tilde{f}_{B}(Q_{\epsilon})}{\epsilon\Gamma}\right)Gdx dt<\epsilon_1,
\end{eqnarray}
then
\begin{eqnarray}
\sup_{P_{\delta R}}\left(\frac{L_1}{2}|\nabla Q_{\epsilon}|^2+\frac{\tilde{f}_{B}(Q_{\epsilon})}{\epsilon}\right)\le C(\delta R)^{-2}
\end{eqnarray}
where the constant $\delta>0$ depends only on $e_{\epsilon}(Q_0,\nabla Q_0)$ and $\min\{R,1\}$.
\end{Lemma}

Proof: We follow the same line as in \cite[Theorem 5.1]{struwe88}.
Let $r_{1}=\delta R$, $\delta\in (0,1/2)$ to be determined later. For $r,\sigma\in (0,r_1)$, $r+\sigma<r_1$, and any $z_0=(x_0,t_0)\in P_{r}$, we have
\beno
I(\sigma,z_0)\doteq\sigma^{-3}\int_{P_{\sigma}(z_0)}e_{\epsilon}(Q_{\epsilon},\nabla Q_{\epsilon})dxdt&\le& C\int_{P_{\sigma}(z_0)}e_{\epsilon}(Q_{\epsilon},\nabla Q_{\epsilon})G_{(x_0,t_0+2\sigma^2)}dxdt\\
&\le &C\int_{T_{\sigma}(t_0+2\sigma^2)}e_{\epsilon}(Q_{\epsilon},\nabla Q_{\epsilon})G_{(x_0,t_0+2\sigma^2)}dxdt.
\eeno
Moreover, apply Remark \ref{rem:monotonicity}, choose $\delta$ small enough, and take  $t_0+2\sigma^2-4R_1^2=-R^2$ and $t_0+2\sigma^2-4R_2^2=-4R^2$, then we deduce that
\beno
I(\sigma,z_0)&\le&C\min\left(\int_{t_0+2\sigma^2-4R_1^2}^{t_0+2\sigma^2-R_1^2},\int_{t_0+2\sigma^2-4R_2^2}^{t_0+2\sigma^2-R_2^2}\right)\left(e_{\epsilon}(Q_{\epsilon},\nabla Q_{\epsilon})G_{(x_0,t_0+2\sigma^2)}\right)dxdt\\
&\le &C\int_{T_{R}}e_{\epsilon}(Q_{\epsilon},\nabla Q_{\epsilon})G_{(x_0,t_0+2\sigma^2)}dxdt.
\eeno
Direct calculation shows that for given $\epsilon_2>0$, if $\delta>0$ is small enough, then we have
\begin{eqnarray*}
G_{(x_0,t_0+2\sigma^2)}(x,t)&\le& C\exp\left(C\delta^2\frac{|x|^2}{4|t|}\right)G(x,t)\nonumber\\
 &\le & \left\{\begin{split}
                 CG(x,t),& ~~~\mbox{if}~~|x|\le \frac{R}{\delta}, \nonumber\\
                  CR^{-3}\exp(-C\delta^{-2}),& ~~~\mbox{if}~~|x|\ge \frac{R}{\delta}
              \end{split}
\right.\nonumber\\
&\le& CG(x,t)+CR^{-2}\exp(-\ln R-C\delta^2)\nonumber\\
&\le &CG(x,t)+\epsilon_2R^{-2},
\end{eqnarray*}
which holds on $T_{R}$, and here $C$ is independent of $\delta$ and $R$. Select $\delta\sim (|\ln R|+|\ln\epsilon_2|)^{-1/2}$ for a small R and independent of $R$ if $R\ge 1$. Thus, it follows that
\begin{eqnarray}\label{eq:e epsilon small}
\sigma^{-3}\int_{P_{\sigma}(z_0)}e_{\epsilon}(Q_{\epsilon},\nabla Q_{\epsilon})dxdt\le C\Psi(R)+C\epsilon_{2} e_{\epsilon}(Q_0,\nabla Q_0)\le C(\epsilon_1+\epsilon_{2} e_{\epsilon}(Q_0,\nabla Q_0)).
\end{eqnarray}

For simplicity, we let $e_{\epsilon}(Q_{\epsilon},\nabla Q_{\epsilon})=e_{\epsilon}(Q_{\epsilon})$. Since $Q_{\epsilon}$ is regular, there exists $\sigma_\epsilon\in (0,r_1)$ such that
\begin{eqnarray}\label{eq:e Q epsilon}
(r_1-\sigma_\epsilon)^2\sup_{P_{\sigma_\epsilon}}e_{\epsilon}(Q_{\epsilon})=\max_{0\le \sigma\le r_1}(r_1-\sigma)^2\sup_{P_{\sigma}}e_{\epsilon}(Q_{\epsilon}).
\end{eqnarray}
Also, there exists a point $(x_\epsilon,t_\epsilon)\in \bar{P}_{\sigma_\epsilon}$  such that
\begin{eqnarray}\label{eq:e epsilon max}
\sup_{P_{\sigma_{\epsilon}}}e_{\epsilon}(Q_{\epsilon})=e_{\epsilon}(Q_{\epsilon})(x_{\epsilon},t_{\epsilon})=e_\epsilon.
\end{eqnarray}
Set $\rho_\epsilon=\frac{1}{2}(r_1-\sigma_\epsilon)$. Then it follows from (\ref{eq:e Q epsilon}) and (\ref{eq:e epsilon max}) that
\beno
\sup_{P_{\rho_\epsilon}(x_\epsilon,t_\epsilon)}e_{\epsilon}(Q_{\epsilon})\le \sup_{P_{\sigma_\epsilon+\rho_\epsilon}}e_{\epsilon}(Q_{\epsilon})\le 4e_\epsilon.
\eeno
Denote
\beno
\tilde{Q}_{\epsilon}(x,t)=Q_{\epsilon}(\frac{x}{\sqrt{e_\epsilon}}+x_\epsilon,\frac{t}{e_\epsilon}+t_\epsilon),
\eeno
which solves the equation (\ref{qtensorflow}) in $P_{r_\epsilon}$ with $\tilde{\epsilon}={\epsilon e_{\epsilon}}$ and $r_{\epsilon}=\sqrt{e_\epsilon}\rho_\epsilon$. Moreover, $\tilde{Q}_{\epsilon}$ satisfies
\begin{eqnarray*}
e_{\tilde{\epsilon}}(\tilde{Q}_{\epsilon})(0,0)=1,\quad
\sup_{P_{r_\epsilon}}e_{\tilde{\epsilon}}(\tilde{Q}_{\epsilon})\le 4.
\end{eqnarray*}
If $r_\epsilon\ge 1$,  $dist(\tilde{Q}_{\epsilon},\mathcal{N})$ convergence uniformly to 0 on $P_{1}$, and there exists $\epsilon_0'$,  such that $dist(\tilde{Q}_{\epsilon},\mathcal{N})<\epsilon_0$ for $\epsilon<\epsilon_0'$. Thus, Lemma \ref{lem:bochner inequlity} implies that
\begin{eqnarray*}
(\partial_{t}-\Delta)e_{\tilde{\epsilon}}(\tilde{Q}_{\epsilon})\le c_1e_{\tilde\epsilon}(\tilde{Q}_{\epsilon}),~~~\mbox{on ~}P_{1}.
\end{eqnarray*}
Moser's Harnack inequality shows that
\beno
1=e_{\tilde{\epsilon}}(\tilde{Q}_{\epsilon})(0,0)\le C \int_{P_1}e_{\tilde{\epsilon}}(\tilde{Q}_{\epsilon})dxdt,
\eeno
while, (\ref{eq:e epsilon small}) tells us
\beno
\int_{P_1}e_{\tilde{\epsilon}}(\tilde{Q}_{\epsilon})dxdt=(\sqrt{e_{\epsilon}})^3\int_{P_{\frac{1}{\sqrt{e_{\epsilon}}}}(x_{\epsilon},t_{\epsilon})}e({Q}_{\epsilon})dxdt\le c(\epsilon_1+\epsilon_{2} e_0(Q_0)),
\eeno
which leads to a contradiction if $\epsilon_1$ and $\epsilon_2$ are suitably small.

Hence, we may assume that $r_\epsilon\le 1$. Then
\beno
1=e_{\tilde{\epsilon}}(\tilde{Q}_{\epsilon})(0,0)&\le& C r_\epsilon^{-5}\int_{P_{r_\epsilon}}e_{\tilde{\epsilon}}(\tilde{Q}_{\epsilon})dxdt=C r_\epsilon^{-2}\rho_\epsilon^{-3}\int_{P_{\rho_\epsilon}}e_{{\epsilon}}({Q}_{\epsilon})dxdt,
\eeno
and using (\ref{eq:e epsilon small}), we get
\beno
\rho_\epsilon^2 e_{\epsilon}=r_\epsilon^{2}\leq  C,
\eeno
then
\beno
\max_{0<\sigma<r_\epsilon}(r_\epsilon-\sigma)^2\sup_{P_\sigma}e(Q_\epsilon)\leq 4\rho_\epsilon^2 e_{\epsilon}\leq C,
\eeno
which implies the required result by choosing $\sigma=\frac12 r_\epsilon=\delta R.$
\endproof

\section{Proof of the main theorem and the equation of $n$}

By Lemma \ref{lem:energy inequality}, we know that for given smooth data $Q_0:\R^3\rightarrow \mathcal{N}$ with $\nabla Q_0\in L^2(\R^3)$, there exist a subsequence of $Q_{\epsilon}$ (also denoted by $Q_{\epsilon}$) and a function $Q(x,t)$, such that as $\epsilon\rightarrow 0$, we have
\ben\label{eq:weak convergence}
&&\nabla Q_{\epsilon}\rightharpoonup\nabla Q~~~~~~\mbox{weakly}^{*}~~\mbox{in}~~ L^{\infty}([0,\infty);L^2(\R^3)),\nonumber\\
&&\partial_tQ_{\epsilon}\rightharpoonup\partial_tQ~~~~~~~\mbox{weakly}~~\mbox{in}~~~L^2(\R^3\times\R_{+}),\nonumber\\
&&Q_{\epsilon}\rightharpoonup Q~~~~~~\mbox{weakly}~~\mbox{in}~~H^{1,2}_{loc}(\R^3\times\R_{+}),\nonumber\\
&&\tilde{f}_{B}(Q_{\epsilon})\rightarrow 0,~~\mbox{in }L^1_{loc}(\R^3\times \R_{+}),
\een
which yield that
\begin{eqnarray}\label{eq:global bound of Q}
\partial_t Q\in L^2(\R^3\times \R_{+}), \nabla Q\in L^{\infty}([0,\infty),L^2(\R^3)),
\end{eqnarray}
and hence also $Q_{\epsilon}\to Q$ a.e. on $\R^3\times\R_{+}$. Also there is a lifting map  $\mn\in \dot{H}^1(\R^3)$ such that
\ben\label{eq:Q lift n}
Q=s_{+}\left(\mn\otimes\mn-\frac{1}{3}Id\right)\in \mathcal{N}.
\een


\subsection{Proof of Theorem \ref{thm:main}: the limit Q-tensor equations}

We follow the standard arguments as in \cite{struwe88} or  \cite{chenstruwe}. Define the singular set by
\beno
\Sigma=\cap_{R>0}\{z_0\in \R^3\times \R_{+}|\liminf_{\epsilon\rightarrow 0}\int_{T_{R}(z_0)}e_{\epsilon}(Q_{\epsilon})dxdt\ge \epsilon_1\}.
\eeno
Then as in \cite[Theorem 6.1]{struwe88}(see also \cite{chenstruwe}), one can show that $\Sigma$ is closed and has locally finite 3-dimensional Hausdorff-measure with respect to the parabolic metric.

For $z_0\notin \Sigma$, there exists a $R_0>0$, and a subsequence of $\epsilon$, which is still denoted by $\epsilon$, such that
\beno
R^{-3}\int_{T_{R_0}(z_0)}e_{\epsilon}(Q_{\epsilon})G_{z_0}dxdt\le \epsilon_1.
\eeno
It follows from Lemma \ref{lem:interior regular} that
\beno
|\nabla Q_{\epsilon}|,\frac{\tilde{f}_{B}(Q_{\epsilon})}{\epsilon}\le C
\eeno
hold uniformly in a uniform neighborhood $\Omega$ of $z_0$. Let $Q$ be the weak limit of (\ref{eq:weak convergence}). Then there exists a subsequence which we denote as $Q_{\epsilon}$ again, such that
\beno
Q_{\epsilon}\to Q,~~~\mbox{in}~~C_{loc}^{0}(\Omega),\label{equi}\\
\nabla Q_{\epsilon}\rightharpoonup \nabla Q,~~~~\mbox{weakly}^{*}~\mbox{in}~L_{loc}^{\infty}(\Omega).
\eeno
Note that ${\tilde{f}_{B}}(Q_{\epsilon})\le C\epsilon$ is a polynomial of $Q_{\epsilon}$,
then the convergence shows $\tilde{f}_{B}(Q)=0$, i.e., $Q\in \mathcal{N}$. Also  $dist^2(Q_{\epsilon},\mathcal{N})\le |Q_{\epsilon}-Q|^2<\epsilon_0$ for $\epsilon$ small enough. Then (\ref{eq:bochner3}) shows that
\begin{eqnarray}
\left(\partial_t-\frac{L_1}{\Gamma}\Delta\right)e_{\epsilon}(Q_{\epsilon})+\frac{1}{\epsilon^2\Gamma^2}\left|\frac{\delta\tilde{f}_B}{\delta Q}\right|^2(Q_{\epsilon})\le C\quad {\rm on }\,\,\Omega,
\end{eqnarray}
which implies that
\beno
\int_{\Omega}\frac{1}{\epsilon^2\Gamma^2}\left|\frac{\delta\tilde{f}_B}{\delta Q}\right|^2(Q_{\epsilon})\phi dxdt\le C(\phi),
\eeno
for all $\phi\in C_{0}^{\infty}(\Omega)$. Moreover, $\frac{1}{\epsilon\Gamma}\frac{\delta\tilde{F}_b}{\delta Q}(Q_{\epsilon})=\frac{1}{\epsilon\Gamma}\mathcal{J}(Q_{\epsilon})$  is uniformly bounded in $L_{loc}^{2}(\Omega)$, also $(\partial_t-\frac{L_1}{\Gamma}\Delta)Q_{\epsilon}$ is uniformly bounded in $L_{loc}^{2}(\Omega)$, and  similar arguments hold for $\partial_t Q_{\epsilon}$ and $\nabla^2 Q_{\epsilon}$. Then we may assume that
\ben\label{eq:jconvergence}
(\partial_t-\frac{L_1}{\Gamma}\Delta)Q_{\epsilon}\to (\partial_t-\frac{L_1}{\Gamma}\Delta)Q,~~~\mbox{weakly in}~~L_{loc}^{2}(\Omega),\nonumber\\
\frac{1}{\epsilon\Gamma}\left|\mathcal{J}(Q_{\epsilon})\right|=\frac{1}{\epsilon\Gamma}\left|\frac{\delta\tilde{f}_b}{\delta Q}(Q_{\epsilon})\right|\to \bar{\lambda},~~~\mbox{weakly in}~~L_{loc}^{2}(\Omega).
\een

The convergence (\ref{equi}) shows that $Q(x,t)\in C(\Omega,\mathcal{N})$, then the lifting $\mn$ in (\ref{eq:Q lift n}) satisfies $\mn\in C(\Omega, S^2)$.
Lemma \ref{lem:Gamma Q in N} tells us that
\beno
\mathcal{J}(Q_{\epsilon})\bot T_{s_{+}(\mn_3^\epsilon\otimes\mn_3^\epsilon-\frac{Id}{3})}\mathcal{N},\quad s_{+}(\mn_3^\epsilon\otimes\mn_3^\epsilon-\frac{Id}{3})=\pi_{\mathcal{N}}(Q_{\epsilon}),
\eeno
 where $\mn_{3}^{\epsilon}$ is the main eigenvector of $Q_{\epsilon}$.

Now we want to prove that the limit $Q$ satisfies the equations of the Q-tensor flow (\ref{eq:limit Q tensor flow}). By the definition of the matrix norm,
 the vector $\mn_{3}^{\epsilon}$ can be estimated as follows:
 \beno
 \|\nabla\mn_{3}^{\epsilon}\|_{L^{\infty}}\le c\|\nabla \pi (Q_{\epsilon})\|_{L^\infty}\le c\|\nabla Q_{\epsilon}\|_{\infty},
 \eeno
As in \cite[Theorem 2]{bz}, we can assume $\mn_3^{\epsilon}\to \mn$  uniformly on $\Omega'\subset\Omega$.

let $\mn_1^\epsilon$ and $\mn_2^\epsilon$ be unit perpendicular vectors in $T_{\mn_3^{\epsilon}}S^2$, which also continuously depend on $Q_{\epsilon}$. Then the following three vectors are the basis of $(T_{\pi_{\mathcal{N}}(Q_{\epsilon})})^\bot:$
\beno
e_1^{\epsilon}(x,t)&=&\frac{1}{\sqrt{2}}\left(\mn_2^{\epsilon}\otimes\mn_1^{\epsilon}+\mn_1^{\epsilon}\otimes\mn_2^{\epsilon}\right),\\
e_{2}^{\epsilon}(x,t)&=&\frac{1}{\sqrt{2}}\left(\mn_{1}^{\epsilon}\otimes\mn_1^{\epsilon}-\mn_{2}^{\epsilon}\otimes\mn_{2}^{\epsilon}\right),\\
e_{3}^{\epsilon}(x,t)&=&\sqrt{6}\left(\frac{1}{2}\mn_1^{\epsilon}\otimes\mn_{1}^{\epsilon}+\frac{1}{2}\mn_{2}^{\epsilon}\otimes\mn_{2}^{\epsilon}-\frac{Id}{3}\right)=-\frac{\sqrt{6}}{2}(\mn_3^{\epsilon}\otimes\mn_{3}^{\epsilon}-\frac{Id}{3}).
\eeno
Then
\beno
\mathcal{J}(Q_{\epsilon})=f_{1}^{\epsilon}(x,t)e_{1}^{\epsilon}(x,t)+f_{2}^{\epsilon}(x,t)e_{2}^{\epsilon}(x,t)+f_{3}^{\epsilon}(x,t)e_{3}^{\epsilon}(x,t).
\eeno
Note that (\ref{eq:jconvergence}) also shows that
\beno
\frac{f_{i}^{\epsilon}(x,t)}{\epsilon\Gamma}\to \bar{\lambda}_1^i(x,t), ~~~\mbox{weakly in }L_{loc}^{2}(\Omega),
\eeno
for $i=1,2,3$. Due to $Q_{\epsilon}\to Q$ in $C_{loc}^{0}(\Omega)$ uniformly, we can assume $\mn_{2}^{\epsilon}\to \mn_2$ and $\mn_{1}^{\epsilon}\to \mn_1$, where $\mn_2$ and $\mn_1$ are perpendicular on $T_{\mn}S^2$. Thus,
\beno
e_i^{\epsilon}(x,t)\to e_{i}(x,t), ~~~\mbox{in~} C_{loc}^{0}(\Omega), ~i=1,2,3,
\eeno
and
\beno
\int_{\Omega'}\frac{\mathcal{J}(Q_{\epsilon})}{\epsilon\Gamma}\cdot\phi dxdt&=&\int_{\Omega'}\frac{f_{1}^{\epsilon}(x,t)}{\epsilon\Gamma}\left(e_{1}^{\epsilon}-e_{1}\right)\cdot\phi dxdt
+\int_{\Omega'}\frac{f_{2}^{\epsilon}(x,t)}{\epsilon\Gamma}\left(e_{2}^{\epsilon}-e_{2}\right)\cdot\phi dxdt\nonumber\\
&&+\int_{\Omega'}\frac{f_{3}^{\epsilon}(x,t)}{\epsilon\Gamma}\left(e_{3}^{\epsilon}-e_{3}\right)\cdot\phi dxdt\to 0
\eeno
for any vector filed $\phi\in L^2_{loc}(\Omega)$ such that $\phi(z)\in T_{Q}\mathcal{N}$ $a.e.$ on $\Omega$, and any $\Omega'\Subset \Omega$, which and (\ref{eq:jconvergence}) yield that
\beno
(\partial_t-\frac{L_1}{\Gamma}\Delta)Q\perp T_{Q}\mathcal{N}, ~~~~a.e.~\mbox{in}~~ \Omega.
\eeno
Then there exists a unit normal vector filed $\gamma_{\mathcal{N}}(Q)\in (T_{Q}\mathcal{N})^{\bot}_{\mathcal{Q}_0}$ along $Q$ and a scalar function $\lambda\in L_{loc}^{2}(\Omega)$ such that
\ben\label{eq:limit Q tensor flow 2}
\partial_t Q-\frac{L_1}{\Gamma}\Delta Q+\lambda\gamma_{\mathcal{N}}(Q)=0\label{ae}
\een
a.e. on $\Omega$ and in the sense of distribution.

Standard arguments show that $Q$ also weakly solves (\ref{eq:limit Q tensor flow 2}) on $\R^3\times \R_{+}$ (for example, see \cite{struwe88}), and we omitted the details.

\subsection{Proof of Theorem \ref{thm:main}:  the harmonic map flow}

Note that $\mn\in \dot{H}^1(\R^3\times \R,S^2)$ is a lifting such that $Q=s_{+}(\mn\otimes\mn-\frac{Id}{3})$ and $|Q|$ is a constant. For
$\phi\in C_{0}^{\infty}(\R^3\times\R_{+},\R^3)$,
\begin{eqnarray}\label{1}
0&=&\left( \int_{\R^3}\int_{\R_{+}}\left((\phi\cdot\mn)\partial_tQ:Q+\nabla Q:\nabla ((\phi\cdot\mn)Q)\right)dxdt+\int_{\R^3}\int_{\R_{+}}(\phi\cdot\mn)\lambda\gamma_{\mathcal{N}}(Q):Qdxdt\right)\nonumber\\
&=&\int_{\R^3}\int_{\R_{+}}|\nabla Q|^2(\phi\cdot\mn)dxdt+\int_{\R^3}\int_{\R_{+}}(\phi\cdot\mn)\lambda\gamma_{\mathcal{N}}(Q):Qdxdt.\end{eqnarray}
We use $\phi^{n}$ to denote $(\mn\cdot\phi)\mn$. Then
\begin{eqnarray}\label{eq:phi normal}
\mn\otimes(\phi-\phi^n)+(\phi-\phi^n)\otimes\mn\in T_{Q}{\cal{N}}.
\end{eqnarray}

On the other hand, by (\ref{eq:limit Q tensor flow 2}), (\ref{eq:phi normal}) and Lemma \ref{lem:AB} we derive 
\beno
0&=&\int_{\R^3}\int_{\R_{+}}\left(\partial_t Q-\Delta Q+\lambda\gamma_{\mathcal{N}}(Q)\right):(\mn\otimes\phi+\phi\otimes\mn)\\
&=&-\int_{\R^3}\int_{\R_{+}}Q:\partial_{t}(\mn\otimes\phi+\phi\otimes\mn)dxdt+\int_{\R^3}\int_{\R_{+}}\nabla Q:\nabla(\mn\otimes\phi+\phi\otimes\mn)dxdt\\
&&+\int_{\R^3}\int_{\R_{+}}\lambda \gamma_{\mathcal{N}}(Q):(\mn\otimes\phi+\phi\otimes\mn)dxdt\\
&=&-2\int_{\R^3}\int_{\R_{+}}Q:\left(\partial_t\mn\otimes\phi+\mn\otimes\partial_t\phi\right)dxdt+2\int_{\R^3}\int_{\R_{+}}\nabla Q:\left(\nabla\mn\otimes\phi+\mn\otimes\nabla\phi\right)dxdt\\
&&+\int_{\R^3}\int_{\R_{+}}\lambda\gamma_{\mathcal{N}}(Q):(\mn\otimes\phi^{n}+\phi^{n}\otimes\mn)dxdt\\
&=&-2s_{+}\int_{\R^3}\int_{\R_{+}}\frac{2}{3}\mn\cdot\partial_t\phi-\frac{1}{3}\partial_t\mn\cdot\phi dxdt
+2s_{+}\int_{\R^3}\int_{\R_{+}}|\nabla\mn|^2\mn\cdot\phi+\nabla\mn\cdot\nabla\phi dxdt\\
&&+2\int_{\R^3}\int_{\R_{+}}\phi\cdot\mn\lambda\gamma_{\mathcal{N}}(Q):\mn\otimes\mn dxdt.\\
&=&-2s_{+}\int_{\R^3}\int_{\R_{+}}\frac{2}{3}\mn\cdot\partial_t\phi-\frac{1}{3}\partial_t\mn\cdot\phi dxdt
+2s_{+}\int_{\R^3}\int_{\R_{+}}|\nabla\mn|^2\mn\cdot\phi+\nabla\mn\cdot\nabla\phi dxdt\\
&&+2\int_{\R^3}\int_{\R_{+}}\phi\cdot\mn\frac{\lambda}{s_{+}}\gamma_{\mathcal{N}}(Q): Q dxdt.\\
&=&2s_{+}\int_{\R^3}\int_{\R_{+}}\partial_t\mn\cdot\phi dxdt
-2s_{+}\int_{\R^3}\int_{\R_{+}}|\nabla\mn|^2\mn\cdot\phi dxdt+2s_{+}\int_{\R^3}\int_{\R_{+}}\nabla\mn\cdot\nabla\phi dxdt,
\eeno
where in the last step we used (\ref{1}).

Thus the proof of the Theorem \ref{thm:main} is completed. \endproof

\bigskip
\noindent {\bf Acknowledgments.}
Part of this work is carried out when the first author is visiting Princeton university.
Meng Wang is partially supported by  NSFC  10931001. Wendong Wang is supported NSFC 11301048 and "the Fundamental Research Funds for the Central Universities".
Zhifei Zhang is partially supported by NSF of China under Grants 11371037 and 11425103.


\begin{thebibliography}{99}
\bibitem{bm} Ball J., Majumdar, {\it Nematic liquid crystals: from Maier-Saupe to a continum theory}, Mol. Cryst. Liq. Cryst., 525(2010), 1-11.

\bibitem{bz} Ball J., Zarnescu A., {\it Orientability and energy minimization in liquid crystal models}, Arch. Rational Mech. Anal., 202(2011), 493-535.

\bibitem{be} Beris A., Edwards B., {\it Thermodynamics of flowing systems with internal microstructure}, Oxford Engrg. Sci. Ser. 36, Oxford, Newtork, 1994.

\bibitem{chenstruwe} Chen Y., Struwe M., {\it Existence and partial regularity results for the heat flow for harmonic maps}, Math. Z., 201(1989), 83-103.

\bibitem{ez} E. W., Zhang P., {\it A molecular kinetic theory of inhomogeneous liquid crystal flow and the small Deborah number limit}, Methods and Applications of Analysis, 13(2006), 181-198.

\bibitem{fcl} Feng, J, Chaubal C., Leal L.,{\it Closure approximations for the Doi theory: Which to use in simulating complex flows of liquid-crystalline polymers?}, Journal of Rheology, 42(1998), 1095-1109.

 \bibitem{fls} Feng, J., Leal L., Sgalari, {\it A theory for flowing nenmatic polymers with orientational distortion}, Journal of Rheology, 44(2000), 1085-1101.

\bibitem{g} P. G. De Gennes, {\it The physics of liquid crystals}, Clarendon Press, Oxford, 1974.
\bibitem{hlwz} Han J., Luo Y., Wang W., Zhang P., Zhang Z.,  {\it From microscopic theory to macroscopic theory: a systematic study on modeling for liquid crystals},
 Arch. Rational Mech. Anal.,  215(2015), 741-809.

\bibitem{HW} Huang T. and Wang C.,  {\it Blow up criterion for nematic liquid crystal flows},  Comm. Partial Differential Equations,  37(2012), 875-884.

\bibitem{kd} Kuzuu N., Doi M., {\it Constitutive equation for nematic liquid crystals under weak velocity gradient derived from a molecular kinetic equation}, Jounal of the Physical Society of Japan, 52(1983), 3486-3494.

 \bibitem{m} Majumdar A., {\it Equilibrium order parameters of liquid crystals in the Landau-De Gennes theory}, European Journal of Applied Mathematics, 21.02(2010),  181-203.

  \bibitem{mn} Mottram N., Newton C., {\it Introduction to Q-tensor theory}, University of Strethclyde, Department of Mathematics, Research Report, 10(2004).

\bibitem{mz} Majumdar A., Zarnescu A., {\it Landau-de Gennes theory of nematic liquid crystals: The Oseen-Frank limit and beyond}, Arch. Ration. Mech. Anal., 196(2010), 227-280.

 \bibitem{pz1} Paicu M., Zarnescu A., {\it Energy dissipation and regularity for a coupled Navier-Stokes and Q-tensor system}, Arch. Ration. Mech. Anal., 203(2012), 45-67.

 \bibitem{pz2} Paicu M., Zarnescu A., {\it Global existence and regularity for the full coupled Navier-Stokes and Q-tensor system}, SIAM J. Math. Anal.,  43(2011),2009-2049.

 \bibitem{qs} Qian T., Sheng P., Generalized hydrodynamic equations for nematic liquid crystals, Phyical Review E, 58(1998), 7475-7485.

\bibitem{keller} J. Rubinstein, P. Sternberg, J. B. Keller, {\it Reaction-diffusion processes and evolution to harmonic maps},  SIAM J. Appl. Math.,  49(1989),1722-1733.

\bibitem{SchoenU} R. M. Schoen, K. Uhlenbeck, {\it A regularity theory for harmonic maps,} J. Differential Geometry, 17(1982), 307-335.

\bibitem{wzz1} Wang W., Zhang P. , Zhang Z. {\it The small Deborah number limit of the Doi-Onsarger equation to the Ericksen Leslie equation}, Communications on Pure and Applied Mathematics, online.

\bibitem{wzz2} Wang W., Zhang P. , Zhang Z., {\it Rigorous derivation from Landau-de Gennes Theorey to Ericksen-Leslie theory},  SIAM J. Math. Anal.,  47(2015), 127-158.

\bibitem{nz} Nguyen L.,  Zarnescu A.,{\it Refined approximation for minimizers of a Landau-de Gennes energy functional}, Calc. Var.,  4(2013), 383-432.

\bibitem{struwe88} Struwe M., On the evolution of harmonic maps in higher dimension, J. Differential Geometry, 28(1988), 485-502.
\end{thebibliography}
\end{document}